\documentclass[a4paper, 11pt]{amsart} 
\usepackage[bookmarks]{hyperref}
\usepackage{graphicx}
\usepackage{psfrag}
\usepackage{amscd}

\usepackage{amsmath, amscd, amsthm, amssymb, mathrsfs}

\DeclareMathOperator{\im}{im}

\DeclareMathOperator{\GL}{GL}

\newcommand{\R}{\mathbb R}
\newcommand{\C}{\mathbb C}
\newcommand{\Z}{\mathbb Z}

\newcommand{\diff}{\text{\rm d}}
\newcommand{\del}{\partial}

\renewcommand{\P}{\mathbb P}

\theoremstyle{plain}
	\newtheorem{theorem}{Theorem}[section]
	\newtheorem{proposition}[theorem]{Proposition}
	\newtheorem{lemma}[theorem]{Lemma}
	\newtheorem{corollary}[theorem]{Corollary}

\theoremstyle{definition}

\theoremstyle{plain}
	\newtheorem*{theorem*}{Theorem}
	\newtheorem*{proposition*}{Proposition}
	\newtheorem*{lemma*}{Lemma}
	\newtheorem*{corollary*}{Corollary}
	\newtheorem*{conjecture*}{Conjecture}
\theoremstyle{definition}
	\newtheorem*{definition*}{Definition}
	\newtheorem*{remark*}{Remark}
	\newtheorem*{remarks*}{Remarks}

\numberwithin{equation}{section}

\begin{document}

\title{Toric anti-self-dual Einstein metrics via complex geometry}
\author{Joel Fine}
\begin{abstract}
Using the twistor correspondence, we give a classification of toric anti-self-dual Einstein metrics: each such metric is essentially determined by an odd holomorphic function. This explains how the Einstein metrics fit into the classification of general toric anti-self-dual metrics given in an earlier paper \cite{donaldson-fine}. The results complement the work of Calderbank--Pedersen \cite{calderbank-pedersen}, who describe where the Einstein metrics appear amongst the Joyce spaces, leading to a different classification. Taking the twistor transform of our result gives a new proof of their theorem.
\end{abstract}

\maketitle

\section{Introduction}
The twistor correspondence between anti-self-dual four-manifolds and certain complex three-folds makes it possible to use techniques of complex geometry to solve problems in Riemannian geometry.  A previous paper \cite{donaldson-fine} exploits this approach to describe the local geometry of an anti-self-dual four-manifold which admits (the germ of) a conformally Killing two-torus action. 

The relevant details are summarised in \S\ref{review of twistors} below. The upshot is that such conformal classes are determined by pairs $(\tau, \phi)$ where:
\begin{itemize}
\item 
$\tau$ is a holomorphic involution of a neighbourhood of $0\in \C$ with $\tau'(0)=-1$;
\item
$\phi$ is a holomorphic $\C^2$-valued function, defined on the same neighbourhood, which is $\tau$-odd ($\phi \circ \tau = -\phi$) and is such that $\phi'(0), \overline{\phi'(0)}$ are $\C$-linearly independent in $\C^2$.
\end{itemize}

This article address the question of when such a conformal class admits an invariant Einstein representative. From the twistor correspondence, this is equivalent to the twistor space admitting an invariant twisted holomorphic contact structure. This places strong restrictions on $(\tau, \phi)$. The main result (Theorem \ref{main theorem}) is that an Einstein representative exists if and only if:
\begin{itemize}
\item
$\tau$ is the ``standard'' involution $\C \to \C$ given by $\tau(z) = -z$;
\item
There is a choice of basis for $\C^2$ such that the two corresponding components of $\phi$ satisfy 
$$
(z^2+1)\phi'_1(z) + i(z^2-1)\phi'_2(z) = 2.
$$
\end{itemize}

Given any odd holomorphic function $\phi_1$, there is a unique odd holomorphic $\phi_2$ which satisfies the above ODE. Thus it follows from this result that toric anti-self-dual Einstein metrics are, locally at least, determined by a single odd holomorphic function. See Theorem \ref{classification} for a more precise statement.

There has already been much written about toric anti-self-dual Einstein metrics. Indeed, Calderbank--Pedersen \cite{calderbank-pedersen} have given a  different classification of the local geometry. Their work relies on that of Joyce \cite{joyce} who classified \emph{surface-orthogonal} toric anti-self-dual four-manifolds---called \emph{Joyce spaces}---a class general enough to contain the Einstein metrics (at least in the non-hyperk\"ahler case). The Calderbank--Pedersen classification is obtained by determining exactly when a Joyce space admits an Einstein representative.

This approach is not explicitly twistorial, however, and it is not clear directly from the work of Calderbank--Pedersen how the Einstein metrics fit into the picture developed in \cite{donaldson-fine}. The point of view taken here has the advantage that it provides a description not only of the metric, but also the associated twistor space. The price paid, however, is that an explicit description of the metric is not so straightforward from this point of view. 

Of course, the two classifications are equivalent; as is explained in \S\ref{metric formulae}, taking the twistor transform of the above ODE for $\phi$ gives a new proof of Calderbank and Pedersen's result.

\subsubsection*{Acknowledgments} 
I would like to thank Simon Donaldson, Michael Singer, Richard Thomas and Dominic Wright for helpful discussions. 

\section{Review of the twistor correspondence}
\label{review of twistors}

\subsection{The twistor correspondence for Einstein metrics}

The Penrose twistor correspondence gives a one-to-one correspondence between germs of conformal classes of anti-self-dual 4-manifolds $M$ and \emph{twistor spaces} $Z$ which are complex three-folds with certain properties: $Z$ admits a fixed-point-free anti-holomorphic involution $\gamma$ and contains a $\gamma$-invariant rational curve $L$ with normal bundle $\mathcal O(1) \oplus \mathcal O(1)$. The deformations of $L$ in $Z$ are called \emph{twistor lines}, the $\gamma$-invariant deformations \emph{real twistor lines}. In the context of Riemannian geometry this theory was first developed in detail by Atiyah--Hitchin--Singer \cite{atiyah-hitchin-singer}. 

Many statements about the Riemannian geometry of an anti-self-dual 4-manifold have holomorphic interpretations on its twistor space. Two examples of this phenomenon which will be used here are the following.

\begin{theorem}[Pontecorvo \cite{pontecorvo}]\label{Pontecorvo}
A K\"ahler representative of an anti-self-dual conformal class corresponds to a holomorphic section of $K^{-1/2}$ which is compatible with $\gamma$ and not identically zero on real twistor lines.
\end{theorem}

\begin{theorem}[Ward \cite{ward}. See also \cite{hitchin2, besse}]
An Einstein representative of an anti-self-dual conformal class corresponds to a holomorphic section $\theta$ of $T^*Z \otimes K^{-1/2}$, whose restriction to real twistor lines is non-zero.
\end{theorem}

As is explained in, for example \cite{besse}, in the Einstein case, the quantity $\diff \theta \wedge \theta$ is a constant function which can be identified with the scalar curvature. Hence the metric is hyperk\"ahler precisely when $\diff \theta \wedge \theta =0$, i.e., when the holomorphic distribution $E = \ker \theta$ is integrable. An anti-self-dual Einstein metric which is \emph{not} hyperk\"ahler has $\diff \theta \wedge \theta \neq 0$, making $\theta$ a twisted contact structure. 

Hyperk\"ahler four-manifolds with Killing fields are well understood thanks to the work of Gibbons--Hawking \cite{gibbons-hawking}. Accordingly, we focus on the non-hyperk\"ahler case here.

\subsection{The toric twistor correspondence} 

A previous article \cite{donaldson-fine} used twis\-tor theory to classify the local geometry of anti-self-dual four-manifolds $M$ which admit two linearly independent, commuting, conformally Killing vector fields $X_1, X_2$. The main points are briefly reviewed here.

The $X_i$ lift to holomorphic vector fields $\widetilde X_i$ on $Z$ where they generate (the germ of) a $\C^2$-action. At least when $M$ is not hypercomplex, this $\C^2$-action is free near a generic real twistor line $L$. It is not, however, transverse to $L$; the orbits are tangential at precisely the two antipodal points corresponding to almost complex structures on $M$ for which $X_1, X_2$ span a complex line.

These two orbits, through the points $J, \gamma(J) \in L$ say, are tangential to first order, so nearby orbits meet $L$ twice. This determines a pair of holomorphic involutions, one defined near $J$ and another defined near $\gamma(J)$; for $J' \in L$ near $J$, the involution $\tau$ is defined by setting $\tau(J')$ to be the other point of $L$ which lies on the orbit through $J'$. The involution defined near $\gamma(J)$ is related to $\tau$ by  $\gamma$ and so carries no extra information.

The action also determines a $\C^2$-valued holomorphic function $\phi$, defined on the domain of $\tau$. For $J'$ near $J$, $\phi(J')$ is the unique element of $\C^2$ which satisfies $\phi(J')\cdot J' = \tau (J')$ where $\cdot$ denotes the $\C^2$-action. The function $\phi$ is $\tau$-odd, i.e., $\phi\circ\tau = -\phi$. 

Let $z$ be a coordinate on $L$ which is $0$ and $\infty$ at $J$ and $\gamma(J)$ respectively, and in which the antipodal map is $\gamma(z) =-\bar{z}^{-1}$. This choice of coordinate is unique up to rotations. The domain of $\tau$ is now a neighbourhood $U \subset \C$ of the origin. The function $\phi \colon U \to \C^2$ also has the property that $\phi'(0)$ and $\overline{\phi'(0)}$ are linearly independent vectors in $\C^2$. (This follows from the fact that the normal bundle of $L$  is $\mathcal O(1) \oplus \mathcal O(1)$. If the vectors are linearly dependent, $L$ has normal bundle $\mathcal O(2)\oplus \mathcal O$.) 

Given such a pair $(\tau, \phi)$ the corresponding twistor space can be explicitly reconstructed; hence $(\tau, \phi)$ completely determines the twistor space---and so the anti-self-dual metric---and all such pairs arise. Taking into account the rotational freedom in the choice of coordinate $z$ gives the following result. Let $\mathcal S$ denote the set of germs, taken at generic points, of toric anti-self-dual conformal classes, modulo conformal equivalence; let $\mathcal T$ denote the set of all pairs $(\tau, \phi)$ as described above, modulo the action of $S^1$ by rotation.

\begin{theorem}[Donaldson--Fine, \cite{donaldson-fine}]
The sets $\mathcal S$ and $\mathcal T$ are in natural one-to-one correspondence.
\end{theorem}

It is also possible to describe when the toric anti-self-dual conformal class admits an invariant K\"ahler representative, or, equivalently, when $K^{-1/2}$ admits invariant holomorphic sections. It turns out that this happens if and only if $\tau$ lies in a certain one-parameter family $\{\tau_c \}$ and that when this happens the space of invariant holomorphic sections of $K^{-1/2}$ has complex dimension two.

We should mention that there are many other approaches to studying anti-self-dual four-manifolds with Killing fields. In the case of a single Killing field, Jones--Tod \cite{jones-tod} considered the induced geometry on the three-dimensional space of integral curves, showing that anti-self-dual four-mani\-folds with a Killing field correspond to three-dimensional Einstein--Weyl spaces carrying a monopole. This correspondence has been particularly fruitful; see, for example, the results in \cite{calderbank-pedersen2,lebrun4,lebrun3}. 

The first work in the toric case was that of Joyce \cite{joyce} which focuses on the case when the orthogonal distribution to the action is integrable. Joyce's construction---which is briefly reviewed in section \ref{joyces construction}---does not use the Jones--Tod correspondence. Applying this correspondence, however, leads to an alternative approach to the toric case; namely one considers separately the quotients by each Killing field and then relates the resulting Einstein--Weyl spaces and monopoles. This idea is explored by Calderbank--Mason \cite{calderbank-mason}. It also features in Calderbank--Pedersen's approach to toric anti-self-dual Einstein metrics \cite{calderbank-pedersen}.

\section{Toric anti-self-dual Einstein metrics}

The main question considered here is to decide when the conformal class $(\tau, \phi)$ admits an invariant Einstein representative with non-zero scalar curvature.

\subsection{Symmetries of anti-self-dual Einstein metrics}

We begin by recalling some well-known material concerning symmetries of anti-self-dual Einstein metrics. As is mentioned above, in the non-hyper\-k\"ahler case, the twistor correspondence describes an Einstein metric via a twisted holomorphic contact structure on $Z$, i.e., a holomorphic, maximally non-integrable, hyperplane distribution $E \subset TZ$. A \emph{contact vector field} on $Z$ is one whose flow preserves $E$. The following standard result in contact geometry describes the holomorphic contact fields. (The non-holomorphic version of this result is described in \cite{arnold}. It was first applied in the context of twistor theory in \cite{lebrun}.)

\begin{proposition}\label{contact symmetries}
Let $E\subset TZ$ be a twisted holomorphic contact distribution on a complex manifold $Z^{2n+1}$. The quotient map $TZ \to TZ/E$ defines an isomorphism between the space of holomorphic contact fields and $H^0(Z, TZ/E)$.
\end{proposition}

\begin{proof}
The quotient $TZ \to TZ/E$ is a holomorphic 1-form $\theta$ with values in $TZ/E$. The first step is to move to a space where $\theta$ is a genuine 1-form. Let $Z' = (TZ/E)^*\setminus Z$, where $Z$ is thought of as the zero section. ($Z'$ is the so-called ``symplectisation'' of $Z$.) The pull back of $TZ/E$ under $\pi \colon Z' \to Z$ is tautologically trivialised and, consequently, $\pi^* \theta$ is a holomorphic 1-form on $Z'$. It is just the restriction of the canonical 1-form of $T^*Z$ under the embedding $Z' \subset T^*Z$ determined by $\theta$. Moreover, $\diff (\pi^* \theta)$ is a non-degenerate holomorphic 2-form on $Z'$, i.e., $\theta$ embeds $Z'$ as a complex symplectic submanifold of $T^*Z$. This is equivalent to the fact that $E$ is maximally non-integrable.

Now, let $X$ be a holomorphic contact field on $Z$. Since it preserves $E$ it lifts to a holomorphic symplectic vector field on $Z'$. The symplectic structure is exact and the lift, and hence $X$, is determined by its Hamiltonian $\theta(X)\in H^0(TZ/E)$ interpreted as a function on $(TM/E)^*$. Notice that the Hamiltonian is linear along the fibres of $Z' \to Z$. Conversely, any holomorphic function on $Z'$ which is fibrewise linear---i.e., a section $s\in H^0(TZ/E)$---determines a holomorphic Hamiltonian vector field on $Z'$ which descends to a holomorphic contact field on $Z$.
\end{proof}

Applying this and the twistor correspondence in the case of anti-self-dual Einstein metrics gives the following well-known result. It appears explicitly in \cite{lebrun} and is also implicit in Lemma 13.36 of the earlier \cite{besse} (see also \cite{calderbank-pedersen}, Proposition 2.1 and its corollary); the link between the holomorphic sections of $K^{-1/2}$ considered here and the solutions of the twistor equation on $M$ considered in \cite{besse, calderbank-pedersen} is provided by \cite{hitchin}. 

\begin{corollary}[See \cite{besse,calderbank-pedersen,lebrun}]\label{Killing fields}
Let $M$ be an anti-self-dual Einstein manifold with non-zero scalar curvature and let $V$ be the space of Killing fields on $M$. There is a canonical isomorphism $V \otimes \C \cong H^0(Z, K^{-1/2})$, which respects the real structures.
\end{corollary}

Combining this \ref{Killing fields} with Pontecorvo's Theorem (\ref{Pontecorvo}) now gives a twistor-theoretic proof of the following fact. This was originally proved by Tod in \cite{tod} where he also outlines the twistor argument given below, which is due to LeBrun. We repeat the proof for completeness.

\begin{corollary}[Tod, \cite{tod}]
Let $M$ be an anti-self-dual Einstein manifold with a non-zero Killing field. Then $M$ is conformally K\"ahler on a dense open subset.
\end{corollary}
\begin{proof}
If $M$ has zero scalar curvature then $M$ is hyperk\"ahler and there is nothing to prove. If $M$ has non-zero scalar curvature, then by Corollary \ref{Killing fields} there is a real non-zero section  $s \in H^0(K^{-1/2})$. By Pontecorvo's Theorem, \ref{Pontecorvo} above, K\"ahler representatives of an anti-self-dual conformal class are given by precisely such sections. The K\"ahler metric on $M$ is defined only at the points corresponding to real twistor lines on which $s$ does not vanish. Since, on restriction to any twistor line $L$, $K^{-1/2}|_L \cong \mathcal O(2)$ (this is part of the standard twistor theory, see \cite{atiyah-hitchin-singer}), the zero divisor of $s$ meets the generic line in two distinct points. Hence the corresponding K\"ahler metric is defined almost everywhere on $M$. 
\end{proof}

As an aside, note that it is not possible to improve this to globally conformally K\"ahler in general. A well-known example of this phenomenon is given by $\overline{\P^2}$ with the Fubini--Study metric. This is anti-self-dual and Einstein with Killing fields, yet, since $b^+ = 0$, cannot be K\"ahler. The K\"ahler metric produced by the above result is defined on the complement of a line. In other words, the K\"ahler metric is defined on the blow-up of $\C^2$ at the origin and, in fact, it can be seen to be the Burns metric. This example is explained in detail in an article of LeBrun \cite{lebrun2}. 

\subsection{From an Einstein metric to the holomorphic conditions}

We now consider the case of a toric anti-self-dual conformal class which admits an invariant Einstein representative, with non-zero scalar curvature. Let $(\tau, \phi)$ be the corresponding pair which defines the twistor space $Z$ with $\C^2$-action. By assumption, $Z$ admits an invariant twisted holomorphic contact structure $\theta \in H^0(T^*Z \otimes K^{-1/2})$. The fact that $\theta$ is invariant places strong restrictions on both $\tau$ and $\phi$ which we now describe. 

There are two important rank-two sub-bundles of $TZ$: the trivial $\underline{\C}^2$ sub-bundle defined by the action and the kernel $E$ of $\theta$. The essential idea is to exploit the $\C^2$-invariance of both the intersection $E \cap \underline\C^2$ and of $E$ itself.

The lifts $\widetilde X_i$ of the Killing fields are linearly independent holomorphic contact fields. It follows from Corollary \ref{Killing fields} that $\xi_i=\theta(\widetilde X_i) \in H^0(K^{-1/2})$ are also linearly independent. They give a family of non-zero sections of $K^{-1/2}$ parametrised up to scale by $\P^1$. This means, in the notation of \S\ref{review of twistors}, that $\tau = \tau_c$ for some $c$. The family of sections includes the $\C^*$ family described in the previous paper \cite{donaldson-fine}, the two missing points of $\P^1$ corresponding to sections each with a double zero on the central twistor line, instead of two distinct simple zeros. 

It follows, in particular, that on the central twistor line $\xi_1, \xi_2$ have no common zeros. If they did, every section in the family would share a zero, which, from the explicit description in \cite{donaldson-fine}, never happens. This means that over $L$, $E \cap \underline\C^2$ is a line sub-bundle of $\underline\C^2$. Since $K^{-1/2}|_L$ has degree 2, the corresponding map $L \to \P^1$ is a double cover. 

\begin{lemma}\label{tau extends}
The involution $\tau$ extends to the whole of the central twistor line $L$.
\end{lemma}

\begin{proof}
$\C^2$-invariance means that $\tau$ must flip the sheets of the double cover $L \to \P^1$ and so extends to the whole of $L$.
\end{proof}

So $\tau$ is the ``standard'' involution, given by $\tau(z) = -z$ in a global coordinate on $L$. There is a natural identification $K^{-1/2}|_L \cong TL$; using this write $\xi_i = q_i \del_z$ for quadratic polynomials $q_i$.

\begin{lemma}
The $q_i$ satisfy $q_i(-z) =q_i(z)$ (i.e., they have no linear term).
\end{lemma}

\begin{proof}
Since the double cover $L \to \P^1$ is $\tau$-invariant, there is a scalar $c \in \C^*$ such that $q_i(-z) = c q_i(z)$. This forces $c^2 = 1$ and, since the $\xi_i$ can't both vanish at $z=0$, we must have $c=1$.
\end{proof}

The next step is to examine the $\C^2$-invariance of $E$. This  forces a relation between the two components of $\phi$.

\begin{lemma}
The components of $\phi$ satisfy the ODE:
$$
q_1(z)\phi'_1(z) + q_2(z)\phi'_2(z)= 2.
$$
\end{lemma}

\begin{proof}
Since $\theta$ is non-zero when restricted to a twistor line, $E$ is always transverse to the central line $L$: over $L$,  $TZ = E \oplus TL$. On the other hand, the $\C^2$-action is also transverse to $L$ except at the points $0, \infty$: over $L-\{0, \infty\}$, $TZ = \underline\C^2 \oplus TL$. Hence, away from $0$ and $\infty$, it is possible to express $E$ as the graph of a map $\chi\colon \underline\C^2 \to TL$. The map $\chi$ is determined by $\theta(\chi(v)) = - \theta(v)$ for all $v \in \C^2$. In fact, on restriction to $L$, $\theta$ is a non-zero section of $T^*L \otimes K^{-1/2}|_L \cong \mathcal O$, i.e., a non-zero constant; for convenience rescale $\theta$ so that this constant is one. Then $\chi = (-q_1\del_z, -q_2\del_z)$ and over $L - \{0, \infty\}$, $E=\im Q$ where $Q \colon
\underline\C^2 \to \underline\C^2\oplus TL$ is given by
$$
Q=
\left(
\begin{array}{cc}
1 & 0 \\
0 & 1\\
- q_1\del_z& -q_2\del_z
\end{array}
\right).
$$

When $z$ is in the domain of $\phi$, acting by $\phi(z)$ gives a map $Z \to Z$ taking $z$ to $-z$; let $D(z) \colon T_zZ \to T_{-z}Z$ denote the derivative of this map. With respect to the splitting $T_zZ = \underline\C^2 \oplus T_zL$, 
$$
D = \left( 
\begin{array}{ccc}
1 & 0 & \diff\phi_1\\
0 & 1 & \diff\phi_2\\
0 & 0 & \diff \tau
\end{array}
\right).
$$
By $\C^2$-invariance, $E_{-z} = D(z)(E_z)$. This means that for $z\neq 0$ the two injections $D(z)Q(z)$ and $Q(-z)$ have the same image. Hence there is a unique isomorphism $A(z)$ of $\C^2$ such that $Q(-z)A(z) = D(z)Q(z)$. Writing this out gives 
$$
A = \left(
\begin{array}{cc}
1 - q_1\phi_1' & -q_2\phi_1'\\
-q_1\phi'_2 & 1 - q_2\phi'_2
\end{array}
\right)
$$
(where we have used $q_i(-z) = q_i(z)$).

Since $D^{-1}(z)=D(-z)$ it follows that $A^{-1}(z)=A(-z)$. Note also that $A(-z) = A(z)$, so $A^2=1$. The off-diagonal entries of $A^2$ are automatically zero, whilst the diagonal terms are one if and only if $q_1\phi_1' + q_2\phi_2' = 2$.
\end{proof}

The quadratic polynomials $q_1, q_2$ have no common zeros and so are a basis for the space of quadratics of the form $az^2 + b$. They are not invariants of $Z$ alone since they depend on the choice of basis $\widetilde X_1, \widetilde X_2$ for the $\C^2$-action. This dependence is $\GL(2,\C)$-equivariant, for the obvious action on bases, so by changing basis for the $\C^2$-action, we can arrange for $q_1, q_2$ to be any basis we like.

If we want to maintain compatibility with the real structure on $Z$, however, we should only change basis for the $\R^2$-action on $M$. In this case the vector fields $q_i \del_z$ remain compatible with the antipodal map $z \mapsto -\bar z^{-1}$. By changing basis for the $\R^2$-action we can still arrange for 
$$
q_1(z) = z^2+1, \quad q_2(z) =i(z^2-1).
$$

We now collect the results of this section together:

\begin{proposition}\label{necessary conditions}
If the conformal class corresponding to $(\tau, \phi)$ admits an invariant Einstein representative then:
\begin{itemize}
\item
There is a global coordinate on $L$ in which $\tau(z) = -z$;
\item
There is a choice of basis for the Killing fields such that the components of $\phi$ satisfy the ODE
$$
(z^2+1)\phi'_1(z) + i(z^2-1)\phi'_2(z) = 2.
$$
\end{itemize}
\end{proposition}

\subsection{From the holomorphic conditions to an Einstein metric}

It remains, of course, to show that the conditions of Proposition \ref{necessary conditions} are sufficient. Let $(\tau, \phi)$ be holomorphic data determining a toric anti-self-dual conformal class where $\tau(z) = -z$ in a global coordinate. As is explained in \cite{donaldson-fine}, the corresponding twistor space $Z$ is built by gluing three pieces together.  Assume $\phi$ is defined on the disc of radius $r <1$. Let $D_1 = \{z : |z| <r\}$, $\Omega = \{ z :  r/2 < |z| < 2/r\}$ and $D_2 = \{z : 1/r<|z|\}$ and let $U$ denote a small neighbourhood of the origin in $\C^2$. Denote by $X_i$ the quotient of $D_i \times U$ under $(z, v) \sim (-z, v + \phi(z))$, and let $X_\Omega = \Omega\times U$. Provided $U$ is chosen small enough, the quotient map is an isomorphism near the boundary of each $D_i \times U$ and this is used to glue the three pieces to form $Z = X_1 \cup X_\Omega \cup X_2$. 

\begin{proposition}\label{sufficient conditions}
Let $q_1,q_2$ be a basis for the space of quadratic polynomials of the form $az^2 + b$. If $\tau(z) =-z$ in a global coordinate on $L$ and 
$$
q_1(z)\phi'_1(z) + q_2(z)\phi'_2(z)= 2.
$$
then the twistor space corresponding to $(\tau, \phi)$ admits an invariant twisted contact structure $\theta \in H^0(T^*Z\otimes K^{-1/2})$ which is non-zero on twistor lines.
\end{proposition}

\begin{proof}
The first step is to define the kernel $E \subset TZ$. Let
$$
V_j= \frac{\del}{\del v_j} + q_i\frac{\del}{\del z}
$$
and consider the hyperplane distribution $\langle V_1, V_2\rangle$ on the parts $D_i \times U$ and $X_\Omega =\Omega\times U$. 

$X_i$ is the quotient of $D_i \times U$ by the $\Z_2$-action generated by $g(z,v) = (-z, v + \phi(z))$. Now
\begin{eqnarray*}
g_*\left(\frac{\del}{\del v_j}\right) 
&=& 
\frac{\del}{\del v_j },\\
g_*\left(\frac{\del}{\del z}\right) 
&=&
-\frac{\del}{\del z} 
+ \phi'_1\frac{\del}{\del v_1}
+ \phi'_2\frac{\del}{\del v_2}
\end{eqnarray*}
Direct calculation gives
$$
g_*(V_1)
=
(1 + q_1)\phi'_1 V_1
+
q_1 \phi'_2V_2,
$$
and so $g_*(V_1) \in \langle V_1, V_2\rangle$. The point is that $V_1, V_2$ and $\del/\del z$ form a basis for the tangent bundle. The hypothesis $q_1\phi_1' + q_2\phi_2' = 2$ implies that, in this basis, the $\del/\del z$-component of $g_*(V_1)$ vanishes. Similarly, $g_*(V_2) \in \langle V_1, V_2\rangle$. This means that the distribution descends to each $X_i$. Under the gluing they fit together with the distribution on $X_\Omega$ to give a rank two sub-bundle $E \subset TZ$. $\C^2$-invariance follows from the invariance of the distribution on $D_i \times U$.

The distribution $E$ defines a $\C^2$-invariant $TZ/E$-valued holomorphic one-form $\theta$. To check that $\theta$ is a contact form, notice that the Lie bracket $[V_1, V_2]$ upstairs has simple zeros at $z=0$ and $z=\infty$ and is non-zero elsewhere. Hence, on passing to the quotients, $X_1, X_\Omega, X_2$, the distribution $E$ is everywhere non-integrable.

The next step is to identify $TZ/E$ with $K^{-1/2}$. The canonical divisor of $Z$ is a pair of $\C^2$-orbits with multiplicity; more precisely, let $p\in L$ with $\tau(p) \neq p$ and let $V_p\subset Z$ denote the $\C^2$-orbit of $p$; then $K = -2(V_p + V_{\tau(p)})$. To see this, denote by $f \colon L \to \P^1$ the branched double covering $f(z) =z^2$ corresponding to $\tau$. The three open sets $D_1, D_2, \Omega$ cover $\P^1$. By moving $p$ if necessary, assume that $f(p) \in \Omega$. Let $\alpha$ be a meromorphic one-form on $\P^1$ with a double pole at $f(p)$ and holomorphic elsewhere. Define a meromorphic three-form $\beta$ on $Z$ by taking $\beta = \alpha\wedge\diff v$ on $X_1$ and $X_2$, and $\beta = f^*\alpha\wedge\diff v$ on $X_\Omega$. It has double poles at $f^{-1}(f(p)) = \{p, \tau(p)\}$ and so the canonical divisor is as claimed.

To show that $TZ/E \cong K^{-1/2}$ it suffices to find a holomorphic section which vanishes to order one along $V_{p}+V_{\tau(p)}$ for some $p$. This is done by projecting, say $\del/\del v_1$. Modulo $E$, this is equal to $q_1\del/\del z$ which has simple zeros along the orbits through the roots of $q_1$, which are $\tau$-invariant.

Finally, notice that, by construction, $E$ is transverse to $L$ and hence, on restriction to $L$, and so also on other nearby twistor lines, $\theta$ is non-zero. 
\end{proof}

To produce a contact structure compatible with the real involution it is necessary to begin with $q_i \del_z$ compatible with $\gamma$. It is then straightforward to check that the distribution $\langle V_1,V_2 ,\rangle$ is real: it is invariant under the real structure on $\Omega \times U$ and taken to itself under the antiholomorphic map $D_1\times U \to D_2 \times U$. This means that the contact structure on $Z$ is also real. Taken together with Proposition \ref{necessary conditions} this gives the following result.

\begin{theorem}\label{main theorem}
The toric anti-self-dual conformal class corresponding to the pair $(\tau, \phi)$ admits an invariant Einstein representative if and only if:
\begin{itemize}
\item There is a global coordinate on $L$ in which $\tau(z) = -z$;
\item There is a choice of basis for the Killing fields such that the components of $\phi$ satisfy the ODE
$$
(z^2+1)\phi'_1(z) + i(z^2-1)\phi'_2(z) = 2.
$$
\end{itemize}
\end{theorem}

This can be rephrased as a ``classification'' of toric anti-self-dual Einstein metrics with non-zero scalar curvature. Let $\mathcal E$ be the set of germs of toric anti-self-dual Einstein metrics taken at generic points, with non-zero scalar curvature. Let $\mathcal F$ denote the set of odd $\C$-valued holomorphic functions $\psi$, defined on a neighbourhood of $0\in\C$ and which have $|\psi'(0)-1| \neq 1$, modulo the action of $S^1$ by rotations in $\C$.

\begin{theorem}\label{classification}
The sets $\mathcal E$ and $\mathcal F$ are in natural one-to-one correspondence.
\end{theorem}

\begin{proof}
Given an odd holomorphic function $\phi_1$ there is a unique odd holomorphic $\phi_2$ for which the ODE is satisfied. (The ODE determines $\phi_2$ up to a constant which is fixed by the requirement $\phi_2(0)= 0$.) In order that this  choice of $\phi$ corresponds to a genuine conformal class it remains to check that $\phi'(0)$ and $\overline{\phi'(0)}$ are linearly independent, which is equivalent to saying that $\overline{\phi'_1(0)}\phi'_2(0)$ is not real.

If $\phi_1, \phi_2$ satisfy the ODE, then $\phi_2'(0) =-i(2-\phi_1'(0))$. Now $\overline{\phi'_1(0)}\phi'_2(0)$ is real if and only if $2\phi_1'(0) - |\phi'_1(0)|^2$ is imaginary which happens if and only if $|\phi'_1(0) -1 | = 1$. Hence the map $(\tau, \phi) \mapsto \phi_1$ gives a bijection $\mathcal E \to \mathcal F$.
\end{proof} 

The two components of $\mathcal F$ given by $|\psi'(0) - 1|>1$ and $|\psi'(0)-1| <1$ correspond to the sign of the scalar curvature in $\mathcal E$. It is interesting to note, perhaps, that as you cross the circle $|\psi'(0) - 1| = 1$, the twistor spaces remain smooth complex manifolds (provided you avoid $\psi'(0)=0$) with twisted holomorphic contact structures; however, the normal bundle to the central line degenerates from $\mathcal O(1) \oplus \mathcal O(1)$ to $\mathcal O(2) \oplus \mathcal O$.

\section{Metric formulae}\label{metric formulae}

In this section we explain how the ODE for $\phi$ gives a relationship between the metric coefficients in two different natural coordinate systems. In one of these coordinate systems, this recovers the main result of Calderbank--Pedersen \cite{calderbank-pedersen}.

\subsection{Twistor coordinates}

Let $(\tau, \phi)$ define a toric anti-self-dual conformal class, with $\tau(z)=-z$ in a global coordinate. As is explained in \cite{donaldson-fine}, the space of twistor lines has natural coordinates $(r,s,v_1,v_2)$ where $v_i \in \C$ are coordinates for the $\C^2$-action and $r,s \in L$ are near $0, \infty$ respectively. We call these \emph{twistor coordinates}. In these coordinates, the conformal structure is given by
$$
\diff r\diff s + 
\frac{
(A_2\diff v_1 - A_1\diff v_2)(B_2\diff v_1 - B_1\diff v_2)}
{(A_2 B_1 - A_1B_2)^2},
$$
where $A_i, B_i$ are functions of $r,s$ defined as follows. Let 
$$
G_i(r,s) = 
- \frac{1}{4\pi i} \int_{C} 
\frac{z\phi_i(z)}{(z^{2}-r)^{1/2}(z^{2}-s)^{1/2}}\diff z,
$$
where the square-root is defined by cutting the plane between $\pm \sqrt r$ and between $\pm \sqrt s$ and the contour $C$ has two components, one around each cut. Then $A_i = \del G_i / \del r$ and $B_i = \del G_i/ \del s$. (Strictly speaking, $A$ and $B$ are only defined up to sign; $A$ behaves like $s^{-1/2}$ near $s=\infty$ and $B$ like $s^{-3/2}$. The products $AB$ in the formula for the metric are well-defined, however. We will gloss over this in what follows.)

The ODE for $\phi$ gives a PDE for the $G_i$ and hence also a relation between the metric coefficients $A_i, B_i$.

\begin{lemma}\label{pde for G}
Suppose that $(az^2 + b)\phi_1'(z) + (cz^2 + d)\phi_2'(z) = 2$. 
Then
$$
(ar + b)\frac{\del G_1}{\del r}
+
(cr + d)\frac{\del G_2}{\del r}
+
(as + b)\frac{\del G_1}{\del s}
+
(cs + d)\frac{\del G_2}{\del s}
=
0.
$$
\end{lemma}

\begin{proof}
This is a direct calculation. For notational convenience, put $\chi = (z^2-r)^{-1/2}(z^2 -s)^{-1/2}$. By assumption, 
$$
2z(a\phi_1 + c\phi_2) \diff z = \diff[(az^2 + b) \phi_1 + (cz^2 + d)\phi_2] - 2\diff z.
$$
Hence,
$$
2aG_1 + 2c G_2 
=
\int_C \diff [(az^2 + b) \phi_1 + (cz^2 + d)\phi_2]\, \chi \diff z,
$$
since $\int_C \chi \diff z = 0$ as $\chi$ is holomorphic between the two components of $C$. Integrating by parts and using $\diff \chi = -(2z^2 - r -s) z \chi^3 \diff z$ gives
$$
2aG_1  + 2cG_2
=
-\frac{1}{4\pi i}\int_C
\frac{[(az^2 + b) \phi_1 + (cz^2 + d)\phi_2] (2z^2 - r -s)}{(z^2-r)(z^2-s)} 
z\chi \diff z.
$$

Taking the $\phi_1$ term first, and expanding in partial fractions, part of the integrand is
$$
\left(
2a 
+
\frac{2(ar + b)}{z^2 - r}
+
\frac{a(s-r)}{z^2 - s}
+\frac{(ar+b)(s-r)}{(z^2-r)(z^2-s)}
\right)
z \phi_1\chi\diff z.
$$
Using $\frac{\del G_1}{\del r} = \frac{1}{2}\int_C (z^2-r)^{-1}z\phi_1\chi\diff z$ and similar formulae for $\frac{\del G_1}{\del s}$ and $\frac{\del^2 G_1}{\del r\del s}$, the $\phi_1$ contribution to the integral is
$$
2a G_1 
+
4(ar + b)\frac{\del G_1}{\del r}
+ 
2(s-r)a \frac{\del G_1}{\del s}
+
4(ar + b)(s-r)\frac{\del^2 G_1}{\del r \del s}.
$$
Next use the second order PDE for $G_1$, 
$$
\frac{\del^2 G_1}{\del r\del s}
=
\frac{1}{2(r-s)}\left(
\frac{\del G_1}{\del r} - \frac{\del G_1}{\del s}
\right),
$$
to write the $\phi_1$ contribution as
$$
2a G_1 
+ 
2(ar + b)\frac{\del G_1}{\del r} 
+ 
2(as + b)\frac{\del G_1}{\del s}.
$$

There is a similar formula for the $\phi_2$ contribution with $a,b$ replaced by $c,d$ and $G_1$ replaced by $G_2$. Putting the pieces together and canceling $2aG_1 + 2c G_2$ from each side gives the result.
\end{proof}

\subsection{Joyce's construction and the Calderbank--Pedersen classification}\label{joyces construction}

A special type of toric anti-self-dual conformal class arises when the orthogonal distribution to the action is integrable; these are called \emph{surface-orthogonal}. It follows from the explicit formula for the conformal class given in the previous section, that when $\tau(z)= -z$ in a global coordinate, any conformal class corresponding to $(\tau, \phi)$ is surface orthogonal (and, indeed, the converse is true). In particular we see, by Lemma \ref{tau extends}, that toric anti-self-dual Einstein metrics are surface orthogonal. This was first proved by Calderbank and Pedersen by different methods (Proposition 3.2 of \cite{calderbank-pedersen}).

In \cite{joyce}, Joyce classifies the local geometry of surface orthogonal anti-self-dual conformal classes. In \cite{calderbank-pedersen} Calderbank and Pedersen give the following description of Joyce's work.

Generically, surface-orthogonal anti-self-dual conformal classes are determined by a pair of axially-symmetric harmonic functions defined on an open set in ${\R}^{3}$. So we take as input a vector-valued function $F=F(x,y)$ on an open set in the upper-half plane $y>0$ satisfying the equation
\begin{equation}\label{ash}
\frac{\del^2F}{\del x ^2}+\frac{\del^2 F}{\del y ^2}+ \frac{1}{y}\frac{\del F}{\del y}= 0.
\end{equation}
Here $x,y$ are real coordinates and $F$ takes values in $\R^{2}$. Set $P= - y \frac{\del F}{\del x}$ and $Q= y \frac{\del F}{\del y}$ so that the derivatives of $P,Q$ satisfy the system of linear equations, the \emph{Joyce equations}:
$$
\frac{\del P}{\del x}= \frac{\del Q}{\del y}\ \ ,\ \  
\frac{\del P}{\del y}+ \frac{\del Q}{\del x}= y^{-1} P. 
$$
Conversely, any solution of these equations arises (at least locally) from some $F$. Now write $P_{i}, Q_{i}$ for the components of $P,Q$. Then Joyce's metric has the form
$$
\frac{\diff x^{2}+\diff y^{2}}{y^{2}} 
+ 
\frac{
(P_{2} \diff u_{1}-P_{1} \diff u_{2})^{2}
+
(Q_{2} \diff u_{1}- Q_{1} \diff u_{2})^{2}}
{(P_{1} Q_{2} - Q_{1} P_{2})^{2}}. $$
This is an anti-self-dual Riemannian metric on a real 4-dimensional manifold; the anti-self-duality condition is a consequence of the Joyce equations.

Using Joyce's classification, Calderbank--Pedersen prove the following result. 

\begin{theorem}[Calderbank--Pedersen, \cite{calderbank-pedersen}]\label{cpe}
A toric anti-self-dual Einstein metric with non-zero scalar curvature is surface orthogonal. Moreover, the two components of $F$ in Joyce's construction satisfy the PDE
$$
x \frac{\del F_1}{\del x}
+ 
y \frac{\del F_1}{\del y}
=
-
\frac{\del F_2}{\del y}.
$$
Conversely, any such pair of axially-symmetric harmonic functions determine a toric anti-self-dual Einstein metric, at least where $P_1Q_2 - Q_1P_2 \neq 0$.
\end{theorem}

Changing from Joyce coordinates to twistor coordinates turns the PDE for $F_1, F_2$ in Theorem \ref{cpe} into the PDE for $G_1, G_2$ in Lemma \ref{pde for G}. To see this we first recall the coordinate transformation, which is described in detail \cite{donaldson-fine}.

We begin by complexifying Joyce's construction to produce an anti-self-dual conformal structure on a \emph{complex} four-manifold, corresponding to the space of twistor lines. To do this simply let $x, y$ be complex and replace $F$ by a $\C^2$-valued function satisfying the same PDE; the conformal class then has same form as above. 

Now let $\zeta = x + iy$, $\xi = x - i y$ and define
$r=(\zeta - i)/(\zeta + i)$, $s=(\xi - i)/(\xi + i)$, to change to twistor coordinates. This transforms Joyce's conformal class into one determined by a function $G$ which is related to $F$ by
$$
\frac{\del F}{\del \zeta}
=
\left(\frac{\zeta + i}{\xi + i}\right)^{1/2}
\frac{\del G}{\del \zeta},
\quad
\frac{\del F}{\del \xi}
=
\left(\frac{\xi + i}{\zeta + i}\right)^{1/2}
\frac{\del G}{\del \xi}.
$$
(An unfortunate misprint in the published version of \cite{donaldson-fine} has the factors $(\zeta+i)^{1/2}(\xi+i)^{-1/2}$ wrongly inverted.)

It is now straightforward to turn Calderbank and Pedersen's PDE in to ours. Changing coordinates from $(x,y)$ to $(\zeta, \xi)$, the equation for $F_1, F_2$ becomes
$$
\zeta\frac{\del F_1}{\del \zeta} 
+ 
\xi\frac{\del F_1}{\del \xi}
+
\frac{\del F_2}{\del \zeta}
+
\frac{\del F_2}{\del \xi}
=
0.
$$
This is equivalent to the following equation for $G_1, G_2$:
$$
(\zeta+i)\left(
\zeta \frac{\del G_1}{\del \zeta}
+ 
\frac{\del G_2}{\del \zeta}
\right)
+
(\xi+i)\left(
\xi \frac{\del G_1}{\del \xi}
+ 
\frac{\del G_2}{\del \xi}
\right)
=
0.
$$
Finally changing coordinates from $(\zeta, \xi)$ to $(r,s)$ this becomes
$$
(r+1)\frac{\del G_1}{\del r}
+
i(r-1)\frac{\del G_2}{\del r}
+
(s+1)\frac{\del G_1}{\del s}
+
i(s-1) \frac{\del G_2}{\del s}
=
0,
$$
which is precisely the PDE of Lemma \ref{pde for G} with the appropriate values of $a,b,c,d$ as chosen in Theorem \ref{main theorem}.

This calculation shows that our twistor analysis leads to a new proof of one direction of Theorem \ref{cpe}, namely that the PDE for $F_1, F_2$ is necessary. To show the converse, we need a converse to Lemma \ref{pde for G}.

\begin{lemma}\label{converse}
Suppose that
$$
(ar + b)\frac{\del G_1}{\del r}
+
(cr + d)\frac{\del G_2}{\del r}
+
(as + b)\frac{\del G_1}{\del s}
+
(cs + d)\frac{\del G_2}{\del s}
=
0.
$$
Then $(az^2 + b)\phi_1'(z) + (cz^2 + d)\phi_2'(z)$ is constant. 
\end{lemma}

\begin{proof}
Put $f(z)=(az^2 + b)\phi_1'(z) + (cz^2 + d)\phi_2'(z)$. Assuming the PDE for $G$, the same calculation as in the proof of Lemma \ref{pde for G} gives that $\int_C f \chi\diff z=0$. As is explained in \cite{donaldson-fine}, there is an interpretation of this condition in terms of the \v{C}ech cohomology of the elliptic curve associated to $(z^2 - r)^{1/2}(z^2-s)^{1/2}$. Namely, it says that $f$ can be written as the difference of two holomorphic functions $f_1$, $f_2$ each defined on opposite branches of the elliptic curve; moreover, the $f_i$ are unique up to the addition of a constant.

In \cite{donaldson-fine}, this was applied to the odd function $\phi$ itself. Here, however, $f$ is even and this means that $f$ can be extended holomorphically over the whole elliptic curve and hence is constant. To see this, define $\hat f_1(z) = - f_2(-z)$, $\hat f_2(z) = -f_1(-z)$. Now, for $z$ in the domain of $f$, $\hat f_2(z) - \hat f_1(z) = f_2(-z) - f_1(-z) = f(-z)= f(z)$. Hence, by uniqueness, $\hat f_i(z) = f_i(z) +c$. This means that $f(z) = f_2(z) - f_1(z) = -f_1(-z)-f_1(z) - c$ and so $f$ extends over one branch; similarly, $f(z)= f_2(z) + f_2(-z) + c$ and so $f$ extends over the other branch.
\end{proof}

It follows from this that if $F_1, F_2$ satisfy the conditions of Theorem \ref{cpe}, so that $G_1, G_2$ satisfy the conditions of 
Lemma \ref{converse}, then $\phi$ satisfies an ODE of the kind required for the corresponding twistor space to admit an invariant twisted holomorphic contact structure. In other words, Lemmas \ref{pde for G} and \ref{converse} combined with Theorem \ref{main theorem} give a twistor theoretic proof of Calderbank and Pedersen's result.

\bibliographystyle{habbrv}
\bibliography{toric_asd_einstein}

{\small \noindent {\tt joel.fine@imperial.ac.uk }} \newline
{Department of Mathematics, Imperial College, London SW7 2AZ. UK.}

\end{document}